\documentclass[12pt]{amsart}
\usepackage{amssymb}

\newcommand{\Bgp}{{\Z^\N}}

\long\def\forget#1\forgotten{}
\newcommand{\issuenumber}{27}
\newcommand{\issuemonth}{April}
\newcommand{\issueyear}{2009}

\newtheorem{issue}{Issue}

\theoremstyle{definition}

\theoremstyle{remark}

\newcommand{\ed}{
\newpage

\section{Unsolved problems from earlier issues}

\begin{issue}
Is $\binom{\Omega}{\Gamma}=\binom{\Omega}{\Tau}$?
\end{issue}

\begin{issue}
Is $\ufin(\cO,\Omega)=\sfin(\Gamma,\Omega)$?
And if not, does $\ufin(\cO,\Gamma)$ imply
$\sfin(\Gamma,\Omega)$?
\end{issue}

\stepcounter{issue}

\begin{issue}
Does $\sone(\Omega,\Tau)$ imply $\ufin(\Gamma,\Gamma)$?
\end{issue}

\begin{issue}
Is $\fp=\fp^*$? (See the definition of $\fp^*$ in that issue.)
\end{issue}

\begin{issue}
Does there exist (in ZFC) an uncountable set satisfying $\sfin(\B,\B)$?
\end{issue}

\stepcounter{issue}

\begin{issue}
Does $X \nin \NON(\cM)$ and $Y\nin\mathsf{D}$ imply that
$X\cup Y\nin \COF(\cM)$?
\end{issue}

\begin{issue}[CH]
Is $\split(\Lambda,\Lambda)$ preserved under finite unions?
\end{issue}

\begin{issue}
Is $\cov(\cM)=\fo$? (See the definition of $\fo$ in that issue.)
\end{issue}

\begin{issue}
Does $\sone(\Gamma,\Gamma)$ always contain an element of cardinality $\fb$?
\end{issue}

\begin{issue}
Could there be a Baire metric space $M$ of weight $\aleph_1$ and a partition
$\mathcal{U}$ of $M$ into $\aleph_1$ meager sets where for each ${\mathcal U}'\subset\mathcal U$,
$\bigcup {\mathcal U}'$ has the Baire property in $M$?
\end{issue}

\stepcounter{issue} 

\begin{issue}
Does there exist (in ZFC) a set of reals $X$ of cardinality $\fd$ such that all
finite powers of $X$ have Menger's property $\sfin(\cO,\cO)$?
\end{issue}

\begin{issue}
Can a Borel non-$\sigma$-compact group be generated by a Hurewicz subspace?
\end{issue}

\begin{issue}[MA]
Is there an uncountable $X\sbst\R$ satisfying $\sone(\BO,\BG)$?
\end{issue}

\begin{issue}[CH]
Is there a totally imperfect $X$ satisfying $\ufin(\cO,\Gamma)$
that can be mapped continuously onto $\Cantor$?
\end{issue}

\begin{issue}[CH]
Is there a Hurewicz $X$ such that $X^2$ is Menger but not Hurewicz?
\end{issue}

\begin{issue}
Does the Pytkeev property of $C_p(X)$ imply that $X$ has Menger's property?
\end{issue}

\begin{issue}
Does every hereditarily Hurewicz space satisfy $\sone(\BG,\BG)$?
\end{issue}

\begin{issue}[CH]
Is there a Rothberger-bounded $G\le\Bgp$ such that $G^2$ is not Menger-bounded?
\end{issue}

\begin{issue}
Let $\cW$ be the van der Waerden ideal.
Are $\cW$-ultrafilters closed under products?
\end{issue}

\begin{issue}
Is the $\delta$-property equivalent to the $\gamma$-property $\binom{\Omega}{\Gamma}$?
\end{issue}

\stepcounter{issue}

\stepcounter{issue}

\general\end{document}}

\newcommand{\Cantor}{{\{0,1\}^\N}}

\newcommand{\fb}{\mathfrak{b}}

\newcommand{\fc}{\mathfrak{c}}
\newcommand{\fd}{\mathfrak{d}}
\newcommand{\fp}{\mathfrak{p}}

\newcommand{\NON}{{\mathsf   {NON}}}
\newcommand{\COF}{{\mathsf   {COF}}}

\newcommand{\cM}{\mathcal{M}}

\newcommand{\op}{\operatorname}
\newcommand{\cov}{\mathsf{cov}}

\newcommand{\R}{\mathbb{R}}

\newcommand{\fo}{\mathfrak{od}}

\renewcommand{\split}{\mathsf{Split}}
\newcommand{\bq}{\begin{quote}}
\newcommand{\eq}{\end{quote}}
\newcommand{\cO}{\mathcal{O}}
\newcommand{\B}{\mathcal{B}}
\newcommand{\BG}{\B_\Gamma}

\newcommand{\BO}{\B_\Omega}

\newcommand{\sone}{\mathsf{S}_1}    \newcommand{\sfin}{\mathsf{S}_\mathrm{fin}}

\newcommand{\ufin}{\mathsf{U}_\mathrm{fin}}

\newcommand{\nin}{\not\in}
\newcommand{\cF}{\mathcal{F}}

\newcommand{\cW}{\mathcal{W}}

\newcommand{\N}{\mathbb{N}}
\newcommand{\Z}{\mathbb{Z}}

\newcommand{\sbst}{\subseteq}
\newcommand{\by}[2]{\par\hfill\emph{#1}, #2}

\newcommand{\Tau}{\mathrm{T}}
\newcommand{\CE}{\textsc{CE}}

\newcommand{\be}{\begin{enumerate}}
\newcommand{\ee}{\end{enumerate}}
\newcommand{\bi}{\begin{itemize}}
\newcommand{\ei}{\end{itemize}}

\newcommand{\general}{\small\vfill\par\noindent\hrulefill\par
\noindent\textbf{Previous issues.} The previous issues of this
bulletin are available online at\\
\texttt{http://front.math.ucdavis.edu/search?\&t=\%22SPM+Bulletin\%22}
\\[0.1cm]
\textbf{Contributions.} Announcements, discussions, and open problems should be emailed
to \texttt{tsaban@math.biu.ac.il}\\[0.1cm]
\textbf{Subscription.}
To receive this bulletin (free) to your e-mailbox, e-mail us.
}

\newcommand{\arXiv}[5]{\subsection{#2}{#4}\par\hfill{\arx{#1}}\par\hfill\emph{#3}}

\newcommand{\AMS}[3]{\subsection{#1}~\par\hfill{\texttt{#3}}\par\hfill\emph{#2}}

\newcommand{\arx}[1]{\texttt{http://arxiv.org/abs/#1}}
\newcommand{\url}[1]{\bq\texttt{#1}\eq}
\newcommand{\online}[1]{The paper is available online at \url{#1}}

\title[$\mathcal{SPM}$ Bulletin \textbf{\issuenumber} (\issuemonth{} \issueyear)]{%
$\mathcal{SPM}$ Bulletin\\[0.5cm]
Issue number \issuenumber: \issuemonth{} \issueyear{} \CE{}}

\begin{document}
\maketitle

\tableofcontents

\section{Editor's note}

A hard-disk (more precisely, disk-on-key) crash I have experienced recently led to loss
of some of the announcements, and some mess in the chronological order of the remaining ones.
Apology for those.

The special issue of Topology and its Applications, dedicated to SPM, did very well in the
downloads statistics: 
Go to \texttt{http://top25.sciencedirect.com/} and choose the journal.

\medskip

Greetings,

\by{Boaz Tsaban}{tsaban@math.biu.ac.il}

\hfill \texttt{http://www.cs.biu.ac.il/\~{}tsaban}

\section{Research announcements}

\arXiv{0812.5033}{Pseudocompact group topologies with no infinite compact subsets}
{Jorge Galindo and Sergio Macario}
{We show that every Abelian group satisfying a mild cardinal inequality admits
a pseudocompact group topology from which all countable subgroups inherit the
maximal totally bounded topology (we say that such a topology satisfies
property $\sharp$). This criterion is used in conjunction with an analysis of
the algebraic structure of pseudocompact groups to obtain, under the
Generalized Continuum Hypothesis (GCH), a characterization of those
pseudocompact groups that admit such a topology.
 We prove in particular that each of the following groups admits a
pseudocompact group topology with property $\sharp$: (a) pseudocompact groups
of cardinality not greater than $2^{2^{\fc}}$; (b) (GCH) connected pseudocompact
groups; (c) (GCH) pseudocompact groups whose torsion-free rank has uncountable
cofinality.
 We also observe that pseudocompact groups with property $\sharp$ contain no
infinite compact subsets and are examples of Pontryagin reflexive precompact
groups that are not compact.}

\arXiv{0901.1688}
{Selective coideals on $(FIN_k^{[\infty]},\leq)$}
{Jos\'e G. Mijares and Jes\'us Nieto}
{A notion of selective coideal on $(FIN_k^{[\infty]},\leq)$ is given. The
natural versions of the local Ramsey property and the abstract Baire property
relative to this context are proven to be equivalent, and it is also shown that
the family of subsets of $FIN_k^{[\infty]}$ having the local Ramsey property
relative to a selective coideal on $(FIN_k^{[\infty]},\leq)$ is closed under
the Souslin operation. Finally, it is proven that such selective coideals
satisfy a sort of canonical partition property, in the sense of Taylor.}

\AMS{Entire functions mapping uncountable dense sets of reals onto each
other monotonically}
{Maxim R. Burke}
{http://www.ams.org/journal-getitem?pii=S0002-9947-09-04924-1}

\AMS{Effective refining of Borel coverings}
{Gabriel Debs; Jean Saint Raymond}
{http://www.ams.org/journal-getitem?pii=S0002-9947-09-04930-7}

\arXiv{0901.3356}{Symmetry and colorings: Some results and open problems}
{T. Banakh, I. V. Protasov}
{We survey some principal results and open problems related to colorings of
algebraic and geometric objects endowed with symmetries.}

\arXiv{0902.0440}
{Many partition relations below density}
{Saharon Shelah}
{We force $2^\lambda$ to be large and for many pairs in the
interval $(\lambda,2^\lambda)$ a stronger version of the polarized
partition relations hold.  We apply this to problems in general
topology.  E.g. consistently, every $2^\lambda$ is successor of
singular and for every Hausdorff regular space $X$, hd$(X) \le
s(X)^{+3}$, {\rm hL}$(X) \le s(X)^{+3}$ and better for $s(X)$ regular, via a
half-graph partition relation.  For the case $s(X) = \aleph_0$ we get
hd$(X)$, hL$(X) \le \aleph_2$ (we can get $\le \aleph_1 < 2^{\aleph_0}$
but in a subsequence work).}

\arXiv{0902.1944}
{Lindelof indestructibility, topological games and selection principles}
{Marion Scheepers and Franklin D. Tall}
{Arhangel'skii proved that if a first countable Hausdorff space is Lindel\"of,
then its cardinality is at most $2^{\aleph_0}$. Such a clean upper bound for
Lindel\"of spaces in the larger class of spaces whose points are ${\sf
G}_{\delta}$ has been more elusive. In this paper we continue the agenda
started in F.D. Tall, On the cardinality of Lindel\"of spaces with points
$G_{\delta}$, Topology and its Applications 63 (1995), 21 - 38, of considering
the cardinality problem for spaces satisfying stronger versions of the
Lindel\"of property. Infinite games and selection principles, especially the
Rothberger property, are essential tools in our investigations.}

\arXiv{0902.2258}
{Locally precompact groups: (Local) realcompactness and connectedness}
{W. W. Comfort and G. Luk\'acs}
{A theorem of A. Weil asserts that a topological group embeds as a (dense)
subgroup of a locally compact group if and only if it contains a non-empty
precompact open set; such groups are called locally precompact. Within the
class of locally precompact groups, the authors classify those groups with the
following topological properties: Dieudonne completeness; local
realcompactness; realcompactness; hereditary realcompactness; connectedness;
local connectedness. They also prove that an abelian locally precompact group
occurs as the quasi-component of a topological group if and only if it is
precompactly generated, that is, it is generated algebraically by a precompact
subset.}

\arXiv{0902.3786}{The group $\op{Aut}(\mu)$ is Roelcke precompact}
{Eli Glasner}
{Following a similar result of Uspenskij on the unitary group of a separable
Hilbert space we show that with respect to the lower (or Roelcke) uniform
structure the Polish group $G= \op{Aut}(\mu)$, of automorphisms of an atomless
standard Borel probability space $(X,\mu)$, is precompact. We identify the
corresponding compactification as the space of Markov operators on $L_2(\mu)$
and deduce that the algebra of right and left uniformly continuous functions,
the algebra of weakly almost periodic functions, and the algebra of Hilbert
functions on $G$, all coincide. Again following Uspenskij we also conclude that
$G$ is totally minimal.}

\newcommand{\intgr}[1]{\left\lfloor{#1}\right\rfloor}
\arXiv{0902.4448}{An infinite combinatorial statement with a poset parameter}
{Pierre Gillibert, Friedrich Wehrung}
{We introduce an extension, indexed by a partially ordered set~$P$ and cardinal numbers~$\kappa$, $\lambda$, denoted by $(\kappa,{<}\lambda)\leadsto P$, of the classical relation $(\kappa,n,\lambda)\rightarrow\nobreak\rho$ in infinite combinatorics. By definition, $(\kappa,n,\lambda)\rightarrow\rho$ holds, if every map $F\colon[\kappa]^n\to[\kappa]^{<\lambda}$ has a $\rho$-element free set. For example, Kuratowski's Free Set Theorem states that $(\kappa,n,\lambda)\rightarrow n+1$ holds if{f} $\kappa\geq\lambda^{+n}$, where $\lambda^{+n}$ denotes the $n$-th cardinal successor of an infinite cardinal~$\lambda$. By using the $(\kappa,{<}\lambda)\leadsto P$ framework, we present a self-contained proof of the first author's result that $(\lambda^{+n},n,\lambda)\rightarrow n+2$, for each infinite cardinal~$\lambda$ and each positive integer~$n$, which solves a problem stated in the 1985 monograph of Erd\H{o}s, Hajnal, M\'at\'e, and Rado. Furthermore, by using an order-dimension estimate established in 1971 by Hajnal and Spencer, we prove the relation $(\lambda^{+(n-1)},r,\lambda)\rightarrow2^{\intgr{\frac{1}{2}(1-2^{-r})^{-n/r}}}$, for every infinite cardinal~$\lambda$ and all positive integers~$n$ and~$r$ with $2\leq r<n$. For example, $(\aleph_{210},4,\aleph_0)\rightarrow 32{,}768$. Other order-dimension estimates yield relations such as $(\aleph_{109},4,\aleph_0)\rightarrow 257$ (using an estimate by F\"uredi and Kahn) and $(\aleph_7,4,\aleph_0)\rightarrow 10$ (using an exact estimate by Dushnik).
}

\arXiv{0903.0659}
{The Schur $\ell_1$ Theorem for filters}
{Antonio Avil\'es, Bernardo Cascales, Vladimir Kadets, Alexander Leonov}
{We study classes of filters $\cF$ on $\N$ such that weak and
strong $\cF$-convergence of sequences in $\ell_1$ coincide. We
study also analogue of $\ell_1$ weak sequential completeness
theorem for filter convergence.}

\arXiv{0904.1567}
{On uniform asymptotic upper density in locally compact abelian groups}
{Szilard Gy. Revesz}
{Starting out from results known for the most classical cases of $\N, \Z^d, \R^d$
or for sigma-finite abelian groups, here we define the notion of asymptotic
uniform upper density in general locally compact abelian groups. Even if a bit
surprising, the new notion proves to be the right extension of the classical
cases of $\Z^d, \R^d$. The new notion is used to extend some analogous results
previously obtained only for classical cases or sigma-finite abelian groups. In
particular, we show the following extension of a well-known result for $\Z$ of
Furstenberg: if in a general locally compact Abelian group G a subset S of G
has positive uniform asymptotic upper density, then S-S is syndetic.}

\AMS{A $c_0$-saturated Banach space with no long unconditional basic sequences}
{J. Lopez-Abad, S. Todorcevic}
{http://www.ams.org/journal-getitem?pii=S0002-9947-09-04858-2}

\arXiv{0903.2374}{Zero subspaces of polynomials on $\ell_1(\Gamma)$}
{Antonio Avil\'es, Stevo Todorcevic}
{We provide two examples of complex homogeneous quadratic polynomials
$P$ on Banach spaces of the form $\ell_1(\Gamma)$. The first
polynomial $P$ has both separable and nonseparable maximal zero
subspaces. The second polynomial $P$
has the property that while the index-set $\Gamma$ is not countable, all zero subspaces of $P$ are separable.}

\arXiv{0904.0816}
{MAD Families and SANE Player}
{Saharon Shelah}
{We throw some light on the question: is there a MAD family
(= a family of infinite subsets of $\Bbb N$, the intersection of any
two is finite) which is completely separable (i.e. any $X \subseteq \Bbb N$ is
included in a finite union of members of the family or includes a
member (and even continuum many members) of the family).
We prove that it is hard to prove the consistency of the negation:
\be
\item If $2^{\aleph_0} < \aleph_\omega$, then there is such a family.
\item If there is no such families then some situation related to pcf holds whose consistency is large; and if ${\frak a} > \aleph_1$ even unknown.
\ee}

\arXiv{0904.0817}
{On the consistency of $\fd_\lambda > \cov_\lambda(\cM)$}
{Saharon Shelah}
{We prove the consistency of: for suitable strongly inaccessible cardinal
lambda the dominating number, i.e.\ the cofinality of $\lambda^\lambda$ is strictly
bigger than $\cov_\lambda(\cM)$, i.e.\ the minimal number of nowhere dense
subsets of $2^\lambda$ needed to cover it. This answers a question of Matet.
}

\arXiv{0904.1389}
{$o$-Boundedness of free topological groups}
{Taras Banakh, Du\v{s}an Repov\v{s}, Lyubomyr Zdomskyy}
{Assuming the absence of $Q$-points (which is consistent with ZFC) we prove that
the free topological group $F(X)$ over a Tychonov space $X$ is $o$-bounded if
and only if every continuous metrizable image $T$ of $X$ satisfies the
selection principle $U_{fin}(O,\Omega)$ (the latter means that for every
sequence $<u_n>_{n\in w}$ of open covers of $T$ there exists a sequence
$<v_n>_{n\in w}$ such that $v_n\in [u_n]^{<w}$ and for every $F\in [X]^{<w}$
there exists $n\in w$ with $F\subset\cup v_n$). This characterization gives a
consistent answer to a problem posed by C. Hernandes, D. Robbie, and M.
Tkachenko in 2000.}



\arXiv{0903.4691}
{Combinatorial and model-theoretical principles related to regularity of
 ultrafilters and compactness of topological spaces. V}
{Paolo Lipparini}
{We generalize to the relations $(\lambda, \mu) \stackrel{\kappa}{\Rightarrow}
(\lambda', \mu')$ and $ \op{alm}(\lambda, \mu) \stackrel{\kappa}{\Rightarrow} \op{alm}(\lambda', \mu')$ some results obtained in Parts II and IV. We also present a multi-cardinal version.}

\AMS{Menger subsets of the Sorgenfrey line}
{Masami Sakai}
{http://www.ams.org/journal-getitem?pii=S0002-9939-09-09887-6}

\arXiv{0904.3104}
{Combinatorial and model-theoretical principles related to regularity of ultrafilters and compactness of topological spaces. VI}
{Paolo Lipparini}
{We discuss the existence of complete accumulation points of sequences in
products of topological spaces. Then we collect and generalize many of the
results proved in Parts I, II and IV.
 The present Part VI is complementary to Part V to the effect that here we
deal, say, with uniformity, complete accumulation points and $ \kappa
$-$(\lambda)$-compactness, rather than with regularity, $[ \lambda, \mu
]$-compactness and $ \kappa $-$ (\lambda, \mu)$-compactness. Of course, if we
restrict ourselves to regular cardinals, Parts V (for $ \lambda = \mu$) and
Part VI essentially coincide.}

\ed